\documentclass[a4paper,10pt,reqno]{amsart}
\usepackage[body={13cm,21cm}]{geometry}
\usepackage[initials]{amsrefs}

\usepackage{amssymb,amscd}
\usepackage[all]{xy}
\usepackage{nicefrac}
\usepackage{units}

\usepackage[mathscr]{eucal}
\usepackage{enumerate}

\numberwithin{equation}{section}
\theoremstyle{plain}
\newtheorem{thm}{Theorem}[section]

\newtheorem{lem}[thm]{Lemma}
\newtheorem{prop}[thm]{Proposition}

\newtheorem{exa}[thm]{Example}
\newtheorem*{problem*}{Problem}
\theoremstyle{definition}
\newtheorem{Def}[thm]{Definition}
\theoremstyle{remark}

\newtheorem{rem}[thm]{Remark}
\newtheorem*{claim*}{Claim}
\newtheorem*{exa*}{Example}
\newtheorem*{rem*}{Remark}
\newtheorem*{rems*}{Remarks}
\newtheorem*{fact*}{Fact}
\newtheorem*{Def*}{Definition}

\def\Q{{\mathbb Q}}
\def\Z{{\mathbb Z}}

\def\G{{\mathbb G}}

\def\ker {{\rm  Ker}}

\DeclareFontFamily{U}{wncy}{}
\DeclareFontShape{U}{wncy}{m}{n}{%
<5>wncyr5%
<6>wncyr6%
<7>wncyr7%
<8>wncyr8%
<9>wncyr9%
<10>wncyr10%
<11>wncyr10%
<12>wncyr6%
<14>wncyr7%
<17>wncyr8%
<20>wncyr10%
<25>wncyr10}{}
\DeclareMathAlphabet{\cyr}{U}{wncy}{m}{n}
\begin{document}

\title[The sum of two integral squares in $\Q(\sqrt{\pm p})$]
{On the sum of two integral squares in quadratic fields
$\Q(\sqrt{\pm p})$}
\author{Dasheng Wei}

\address{ Academy of Mathematics and System Science,  CAS, Beijing
100190, P.R.China}

\email{dshwei@amss.ac.cn}

\date{\today}

\maketitle

 \bigskip
\section*{\it Abstract}

We propose a method for determining which integers can be written as
a sum of two integral squares for quadratic fields $\Q(\sqrt{\pm
p})$, where $p$ is a prime.

\bigskip

{\it MSC classification} : 11D09; 11E12


\bigskip

{\it Keywords} : integral points, ring class field, reciprocity law.

\section*{Introduction} \label{sec.notation}
Gauss first determined which integers can be written as a sum of two
integral squares. And Niven determined  which integers can be
written as a sum of two integral squares for the imaginary quadratic
field $\Q({\sqrt{-1}})$ in \cite{Niv}. Nagell further studied the
question for the twenty quadratic fields $\Q(\sqrt{d})$ in
\cite{Nag1} and \cite{Nag2}, where
$$d=\pm 2,\pm 3,\pm 5,\pm 7,\pm 11,\pm 13,\pm 19,\pm 43,\pm 67,\pm
163.$$ His method essentially depends on the fact that the class
number of $\Q(\sqrt{d},\sqrt{-d})$ is 1 when $d$ is one of the above
integers. However, this method can not apply for general quadratic
fields. Recently, Harari showed that the Brauer-Manin obstruction is
the only obstruction for existence of the integral points of a
scheme over a ring of integers of a number field whose generic fiber
is a principal homogeneous space of tori in \cite{Ha08}. However,
the Brauer-Manin obstruction of tori given in \cite{Ha08} is not
constructive and one can not use that result to determine the
existence of integral points for the scheme. Fei Xu and the author
gave another proof of the result in \cite{WX} and \cite{WX2} which
is constructive. In this paper we apply the method in \cite{WX} for
the question for quadratic fields $\Q(\sqrt{\pm p})$, where $p$ is a
prime.

Notation and terminology are standard if not explained. Let $F$ be a
number field, $\frak o_F$ the ring of integers of $F$, $\Omega_F$
the set of all primes in $F$ and $\infty$ the set of all infinite
primes in $F$. For simplicity, we write $\frak p<\infty$ for $\frak
p\in \Omega_F\setminus \infty$. Let $F_\frak p$ be the completion of
$F$ at $\frak p$ and $\frak o_{F_\frak p}$ be the local completion
of $\frak o_F$ at $\frak p$ for each $\frak p\in \Omega_F$. Write
$\frak o_{F_\frak p}=F_\frak p$ for $\frak p\in \infty$. We also
denote the adele ring (resp. the idele ring) of $F$ by $\Bbb A_F$
(resp. $\Bbb I_F$) and
$$F_{\infty}=\prod_{\frak p\in \infty}F_\frak p .$$

Let $E=F(\sqrt{-1})$ and let $T$ be the torus $R^1_{E/F}(\Bbb
G_m)=\ker[R_{E/F}(\Bbb G_{m,E})\rightarrow \Bbb G_{m,F}]$, here $R$
denotes the Weil's restriction (see \cite{Milne98}). Denote
$\lambda$ to be the embedding from $T$ to $ R_{E/F}(\Bbb G_{m,E})$.
Obviously $\lambda$ induces a natural group homomorphism
$$\lambda_E:\ \ \ T(\Bbb A_F)\rightarrow \Bbb I_E.$$

Let $\bold X_\alpha$ be the affine scheme over $\frak o_{F}$ defined
by the equation $x^2+y^2=\alpha$ for a non-zero integer $\alpha \in
\frak o_F$. The generic fiber of $\bold X_\alpha$ is a principle
homogenous space of the torus $T$. The equation $x^2+y^2=\alpha$ is
solvable over $\frak o_F$ if and only if $\bold X_\alpha(\frak
o_F)\neq \emptyset$.
\begin{Def}
Let $K/E$ be a finite abelian extension. Let $\psi_{K/E}: \Bbb
I_E\rightarrow Gal(K/E)$ be the Artin map. We say that $\alpha$
satisfies the Artin condition of $K$ if there is
$$\prod_{\frak p\leq \infty}(x_\frak p,y_\frak p)\in \prod_{\frak p\leq \infty}
\bold X_\alpha(\frak o_{F_\frak p})\text{  such that  }
\psi_{K/E}(f_E[\prod_{\frak p\leq \infty}(x_\frak p,y_\frak p)])=1
$$ where $1$ is the
identity element of $Gal(K/E)$ and $f_E: \prod_{\frak p\leq \infty}
\bold X_\alpha(\frak o_{F_\frak p})\rightarrow \Bbb I_E$ is defined
by
$$f_E[(x_\frak p,y_\frak p)]= \begin{cases} (x_\frak p+y_\frak p \sqrt{-1}, x_\frak p-y_\frak p\sqrt{-1}) \ \ \
& \text{if $\frak p$ splits in $E/F$} \\
x_\frak p+y_\frak p\sqrt{-1} \ \ \ & \text{otherwise.}
\end{cases} $$
\end{Def}
By the class field theory, it is a necessary condition for $\bold
X_\alpha(\frak o_F) \neq \emptyset$ that $\alpha$ satisfies the
Artin condition of $K$. In fact there is a finite abelian extension
$K/E$ that is independent on $\alpha$, such that the Artin condition
of $K$ is also sufficient for $\bold X_\alpha(\frak o_F) \neq
\emptyset$ (see \cite{WX}).

Let $\bold T$ be the group scheme over $\frak o_F$  defined by
$x^2+y^2=1$, which is an integral model of $T$. Since $\bold T$ is
separated over $\frak o_F$, we can view $\bold T(\frak o_{F_\frak
p})$ as a subgroup of $T(F_{\frak p})$. Furthermore, the following
result is proved in \cite{WX}.
\begin{prop} \label{multiple} Let $K/E$ be a finite abelian extension
 such that the group homomorphism $\widetilde{\lambda}_E$ induced by
$\lambda_E$
$$\widetilde{\lambda}_E: \nicefrac{T(\Bbb A_F)}{T(F)\prod_{\frak p\leq \infty}\bold T(\frak o_{F_\frak p})} \longrightarrow
\nicefrac{\Bbb I_{E}}{E^*N_{K/E}(\Bbb I_K)}
$$ is well-defined and injective, where well-defined means
$$ \lambda_E(T(F)\prod_{\frak p\leq \infty}\bold T(\frak o_{F_\frak p}))\subset E^*N_{K/E}(\Bbb I_K).$$
Then $\bold X_\alpha(\frak o_F)\neq \emptyset$ if and only if
$\alpha$ satisfies the Artin condition of $K$.

\end{prop}

In  this paper, we mainly prove the following result.
\begin{thm} Let $p$ be a prime number and $F$ the quadratic field
$\Q(\sqrt{p})$ or $\Q(\sqrt{-p})$. Then the diophantine equation
$x^2+y^2=\alpha$ is solvable over $\frak o_F$ if and only if
$\alpha$ satisfies the Artin condition of $H_L$, where $H_L$ is the
ring class field corresponding to the order $L=\frak o_{F}+ \frak
o_{F} \sqrt{-1}$.
\end{thm}

\section{The sum of two squares in imaginary quadratic fields }

Let $d$ be a square-free positive integer here $d\geq 2$. Let
$F=\Q(\sqrt{-d})$, $\frak o_{F}$ be the integral ring of $F$ and
$E=F(\sqrt{-1})$. One takes the order $L= \frak o_{F}+ \frak o_{F}
\sqrt{-1}$ inside $E$. Let $H_L$ be the ring class field
corresponding to the order $L$.

\begin{prop} \label{ima} Suppose one of
the following conditions holds:

(1) The equation $x^2+y^2=-1$ has an integer solution in $\frak
o_F$.

(2) The equation $x^2+y^2=-1$ has no local integral solutions at a
place of $F$.

Then the diophantine equation $x^2+y^2=\alpha$ is solvable over
$\frak o_F$ if and only if $\alpha$ satisfies the Artin condition of
$H_L$.
\end{prop}
\begin{proof} (1) First we assume $d\neq 3$.
Let $\frak p$ be a place of $F$ and $L_\frak p$ be the $\frak
p$-adic completion of $L$ inside $E_\frak p = E\otimes_F F_\frak p$.
Recall $T=R^1_{E/F}(\G_{m,F})$ and $\bold T$ is the scheme defined
by the equation $x^2+y^2=1$, we have
$$T(F)=\{\beta \in E^*: \ N_{E/F}(\beta)=1\}$$ and
$$\bold T(\frak o_{F_\frak p})=\{\beta \in L_\frak p^\times: \ N_{E_\frak p/F_\frak p}(\beta)=1\}.$$

Since the ring class field $H_L$ of the order $L$ corresponds to the
open subgroup $E^* ( \prod_{\frak p\leq \infty} L_\frak p^\times)$
of $\Bbb I_E$ by the class field theory, the natural group
homomorphism
$$\widetilde{\lambda}_E : \ \ \ \nicefrac{T(\Bbb A_F)}{T(F)\prod_{\frak p\leq \infty}\bold T(\frak o_{F_\frak p})}
 \longrightarrow \nicefrac{\Bbb I_E}{E^*  \prod_{\frak p\leq \infty} L_p^\times}
$$ is well-defined. By Proposition \ref{multiple}, we only need to
show $\widetilde{\lambda}_E$ is injective.

Suppose there are $$\beta\in E^*  \ \ \ \text{ and} \ \ \ i\in
E_\infty^* \prod_{\frak p<\infty} L_\frak p^\times$$ such that
$\beta \cdot i\in T(\Bbb A_E)$.  Then
$$N_{E/F} (\beta i)=N_{E/F}(\beta) N_{E/F}(i) =1 $$  and $$N_{E/F}(\beta)\in
F^*\cap \prod_{\frak p<\infty} \frak o_{F_\frak p}^\times =\{\pm 1\}
.
$$ If $N_{E/F}(\beta)=1$, one concludes that
$$ N_{E/F}(\beta) =N_{E/F}(i) =1 \ \ \ \Rightarrow \ \ \ \beta\in
T(E) \ \ \ \text{and} \ \ \ i\in \prod_{\frak p\leq \infty}\bold
T(\frak o_{F_\frak p}).$$ So $\beta i\in T(E)\prod_{\frak p\leq
\infty}\bold T(\frak o_{F_\frak p})$.

If $N_{E/F}(\beta)\neq 1$, then $N_{E/F}(\beta)=N_{E/F}(i)=-1$. That
is to say that $x^2+y^2=-1$ has local integral solutions at every
local place of $F$. By the assumption, we have $x^2+y^2=-1$ has an
integral solution $(x_0,y_0)$. Let $$\zeta=x_0+y_0\sqrt{-1}\ \text{
and
 }\ \gamma=\beta \zeta,j=i/\zeta.$$ Then $\beta i=\gamma j$ and
$$ N_{E/F}(\gamma) =N_{E/F}(j) =1 \ \ \ \Rightarrow \ \ \ \gamma\in
T(F) \ \ \ \text{and} \ \ \ j\in \prod_{\frak p}\bold T(\frak
o_{F_\frak p}).$$  So $\beta i=\gamma j\in T(F)\prod_{\frak p\leq
\infty}\bold T(\frak o_{F_\frak p})$. Therefore
$\widetilde{\lambda}_E$ is injective.

(2) If $d=3$, then $\frak o_{F}^\times=<\pm 1,\zeta_3>$, where
$\zeta_3$ is a primitive 3-rd root of unity. Since
$\zeta_3=\zeta_3^4$ is a square, we can give a proof for this case
with similar arguments as above.
\end{proof}

In the rest of this section we consider the case that $d$ is a
prime. First we need the following result that can be found in
\cite{Yo}.
\begin{prop} \label{Yo} Let $p$ be a prime. Then

(1) If $p \equiv 1 \mod 4$, then  $x^2-py^2=-1$ is solvable over
$\Z$.

(2) If $p \equiv -1 \mod 8$, then $x^2-py^2=2$ is solvable over
$\Z$.

(2) If $p \equiv 3 \mod 8$, then $x^2-py^2=-2$ is solvable over
$\Z$.
\end{prop}

Now we can prove the following lemma.
\begin{lem}\label{-p} Let $p$ be a prime and $F=\Q(\sqrt{-p})$. Then

(1) If $p\equiv -1 \mod 8$, then  $x^2+y^2=-1$ is not solvable over
$\frak o_{F_\frak p}$, where $\frak p\mid 2$.

(2) If $p\not\equiv -1 \mod 8$, then $x^2+y^2=-1$ is solvable over
$\frak o_F$.
\end{lem}

\begin{proof} (1) If $p\equiv -1 \mod 8$, then $2$ splits into
$\frak p_1$ and $\frak p_2$ in the field $F/\Q$. So the Hilbert
symbol $$(-1,-1)_{\frak p_1}=(-1,-1)_{\frak p_2}=-1.$$ Therefore the
equation $x^2+y^2=-1$ is not solvable over $\frak o_{F_{\frak p_1}}$
and $\frak o_{F_{\frak p_2}}$.

If $p\equiv 1 \mod 4$ or $p=2$, then $x^2-py^2=-1$ has an integral
solution in $\Z$ by Proposition \ref{Yo}. Choose one solution
$(x_0,y_0)$, we have $x_0^2+(y_0\sqrt{-p})^2=-1$.

If $p\equiv 3 \mod 8$, then $x^2-py^2=-2$ has an integral solution
in $\Z$ by Proposition \ref{Yo}. We can choose one solution
$(x_0,y_0)$ and it's easy to see $x_0,y_0\equiv 1 \mod 2$. So
$$\frac{x_0 \pm y_0\sqrt{-p}}{2} \in \frak o_{F} \text{ and }
(\frac{x_0+y_0\sqrt{-p}}{2})^2+(\frac{x_0-y_0\sqrt{-p}}{2})^2=-1.$$
\end{proof}

By Proposition \ref{ima} and Lemma \ref{-p}, we can now prove the
following result.
\begin{thm} \label{l} Let $p$ be a prime number and
$F=\Q(\sqrt{-p})$. Let $H_L$ be the ring class field corresponding
to the order $L=\frak o_F+\frak o_F\sqrt{-1}$. Then the diophantine
equation $x^2+y^2=\alpha$ is solvable over $\frak o_F$ if and only
if $\alpha$ satisfies the Artin condition of $H_L$.
\end{thm}

\begin{rem}
It is possible that the family of equations $x^2+y^2=\alpha,
\alpha\in \frak o_F \text{ and } \alpha \neq 0$ satisfies Hasse
principle even if the ring class field $H_L$ is not trivial. For
example, the equation $x^2+y^2=\alpha$ satisfies Hasse principle and
$H_L$ is not trivial for $F=\Q(\sqrt{-p})$ with $p=23,31,47,59,71$.
The reason is $H_L=E H$ for the above $p$, where $E=F(\sqrt{-1})$
and $H$ is the Hilbert class field of $F$. If $x^2+y^2=\alpha$ has
local solutions for every place, then $\alpha$ automatically
satisfies the Artin condition of $H_L$ by the class field theory.
\end{rem}

Now we use Theorem \ref{l} to give an explicit example. Let
$F=\Q(\sqrt{-79})$. We write
$N_{F/\Q}(\alpha)=2^{s_1}79^{s_2}p_1^{e_1}\cdots p_g^{e_g}$ for any
$\alpha\in \frak o_F$. Let $D(n)=\{p_1, \cdots, p_g \}$ and
$h(x)=x^3-307x+1772$. Denote
$$\aligned& D_1=\{p\in D(n): (\frac{79}{p})=(\frac{-1}{p})=1 \text{ and }
h(x) \equiv 0 \mod p \text{ is not solvable} \}\cr & D_2= \{p\in
D(n): (\frac{79}{p})=-(\frac{-1}{p})=1 \text{ and } h(x)\equiv 0
\mod p \text{ is not solvable}\}.
\endaligned $$
It's easy to see that $e_i$ is even when $p_i\in D_2$.
\begin{exa} Let $F=\Q(\sqrt{-79})$ and let $\alpha$ be an integer in $F$ with
the above notation. Then $x^2+y^2=\alpha$ is solvable over $\frak
o_F$  if and only if

(1) The equation $x^2+y^2=\alpha$ has integral solutions at every
place of $F$.

(2) The sum $$\sum_{p_i\in D_1}e_i+\sum_{p_i\in D_2}e_i/2\neq 1.$$
\end{exa}


\section{The sum of two squares in real quadratic fields}

Let $d$ be a square-free positive integer and $F=\Q(\sqrt{d})$. Let
$\frak o_{F}$ be the  ring of integers of $F$,  $\varepsilon$ the
fundamental unit of $\frak o_{F}$ and $\varepsilon=a+b\sqrt{d}\text{
with }a,b>0$. Let $E=F(\sqrt{-1})$. One takes the order $L= \frak
o_{F}+ \frak o_{F} \sqrt{-1}$ inside $E$. Let $H_L$ be the ring
class field corresponding to the order $L$.

\begin{prop} \label{real}Suppose one of
the following conditions holds:

(1) The equation $x^2+y^2=\varepsilon$ has an integer solution in
$\frak o_F$.

(2) The equation $x^2+y^2=\varepsilon$ has no local integral
solutions at a place of $F$.

Then the diophantine equation $x^2+y^2=\alpha$ is solvable over
$\frak o_F$ if and only if $\alpha$ satisfies the Artin condition of
$H_L$.
\end{prop}
\begin{proof}
Let $\frak p$ be a place of $F$ and $L_\frak p$ be the $\frak
p$-adic completion of $L$ inside $E_\frak p = E\otimes_F F_\frak p$.
Since the ring class field $K_L$ of the order $L$ corresponds to the
open subgroup $E^* ( \prod_{\frak p\leq \infty} L_\frak p^\times)$
of $\Bbb I_E$ by the class field theory, the natural group
homomorphism
$$\widetilde{\lambda}_E : \ \ \ \nicefrac{T(\Bbb A_F)}{T(F)\prod_{\frak p\leq \infty}\bold T(\frak o_{F_\frak p})}
 \longrightarrow \nicefrac{\Bbb I_E}{E^*  \prod_{\frak p\leq \infty} L_p^\times}
$$ is well-defined. By Proposition \ref{multiple}, we only need to
show $\widetilde{\lambda}_E$ is injective.

Suppose there are $$\beta\in E^*  \ \ \ \text{ and} \ \ \ i\in
E_\infty^* \prod_{\frak p<\infty} L_\frak p^\times$$ such that
$\beta \cdot i\in T(\Bbb A_E)$.  Then
$$N_{E/F} (\beta i)=N_{E/F}(\beta) N_{E/F}(i) =1 $$  and $$N_{E/F}(\beta)\in
F^*\cap \prod_{\frak p<\infty} \frak o_{F_\frak p}^\times =\{\pm
\varepsilon^n\}.
$$ Since $N_{E/F}(\beta)$ is totally positive, we have
$N_{E/F}(\beta)=\varepsilon^n$.

When $n$ is even, let $\gamma=\beta
\varepsilon^{n/2},j=i\varepsilon^{-n/2}$. Then $\beta i=\gamma j$
and
$$ N_{E/F}(\gamma) =N_{E/F}(j) =1 \ \ \ \Rightarrow \ \ \ \gamma\in
T(F) \ \ \ \text{and} \ \ \ j\in \prod_{\frak p}\bold T(\frak
o_{F_\frak p}).$$  So $\beta i=\gamma j\in T(F)\prod_{\frak p\leq
\infty}\bold T(\frak o_{F_\frak p})$.

When $n$ is odd, we have $N_{E/F}(i)=\varepsilon^{-n}$. That is to
say that $x^2+y^2=\varepsilon^{-n}$ has integral solutions at every
local place of $F$. Since $n$ is odd and $\varepsilon\in \frak
o_F^\times$, we have $x^2+y^2=\varepsilon$ has integral solutions at
every place of $F$. By the assumption, we have $x^2+y^2=\varepsilon$
has an integer solution $(x_0,y_0)$. Let $\zeta=x_0+y_0\sqrt{-1}$
and $\gamma=\beta
\varepsilon^{(n-1)/2}\zeta,j=i\varepsilon^{(1-n)/2}\zeta^{-1}$. Then
$\beta i=\gamma j.$ And
$$ N_{E/F}(\gamma) =N_{E/F}(j) =1 \ \ \ \Rightarrow \ \ \ \gamma\in
T(F) \ \ \ \text{and} \ \ \ j\in \prod_{\frak p\leq \infty}\bold
T(\frak o_{F_\frak p}).$$  So $\beta i=\gamma j\in T(F)\prod_{\frak
p\leq \infty}\bold T(\frak o_{F_\frak p})$. Therefore
$\widetilde{\lambda}_E$ is injective.
\end{proof}

In the following, we consider the case that $d$ is a prime number.
\begin{lem} \label{p}Let $p$ be a prime and $F=\Q(\sqrt{p})$. Let $\varepsilon$
be the fundamental unit of $\frak o_{F}$ and
$\varepsilon=a+b\sqrt{p}\text{ with }a,b>0$. Then there is a place
$\frak p$ of $F$, such that the equation $x^2+y^2=\varepsilon$ is
not solvable over $\frak o_{F_\frak p}$.
\end{lem}

\begin{proof}
If $p\equiv 1 \mod 4$ or $p=2$, then $x^2-py^2=-1$ has integral
solutions by Proposition \ref{Yo}. Therefore
$N_{F/\Q}(\varepsilon)=-1$. There exists a real place $\frak p$ of
$F$ such that $|\varepsilon|_\frak p<0$. So the equation
$x^2+y^2=\varepsilon$ is not solvable at the real place $\frak p$.

If $p\equiv 3 \mod 4$, then $x^2-py^2=-1$ is not solvable over $\Z$
by Proposition \ref{Yo}. Therefore $N_{F/Q}(\varepsilon)=1$ and
$\varepsilon$ is totally positive. And we know one of the equations
$x^2-py^2=\pm 2$ has an integral solution in $\Z$ by Proposition
\ref{Yo}. We can choose one solution $(x_0,y_0)$ and it is easy to
see  that $x_0$ and $y_0$ are odd. Let
$$A=(x_0^2+py_0^2)/2 \text{ and } B=x_0y_0.$$
Since $x_0,y_0$ are odd, we can see $A,B$ are integers and $B$ is
odd. And $$A^2-pB^2=(x_0^2-py_0^2)^2/4=1.$$ Let
$\varepsilon_1=A+B\sqrt{p}$. Obviously $\varepsilon_1$ is totally
positive and $\varepsilon_1=\varepsilon^m, m\in \Z$.

Let $\frak p$ be the unique place of $F$ over $2$. We assume the
equation $x^2+y^2=\varepsilon$ is  solvable over $\frak o_{F_\frak
p}$. Since $\varepsilon_1=\varepsilon^m$, the equation
$x^2+y^2=\varepsilon_1$ is also solvable over $\frak o_{F_\frak p}$.
For any solution $(x_1,y_1)=(a_1+b_1\sqrt{p}, a_2+b_2\sqrt{p})$, we
have $$(a_1+b_1\sqrt{p})^2+(a_2+b_2\sqrt{p})^2=A+B\sqrt{p}.$$ Then
$$2a_1b_1+2a_2b_2=B.$$
However, we know $B$ is odd, a contradiction is derived.
\end{proof}

By Proposition \ref{real} and Lemma \ref{p}, we can now prove the
following result.
\begin{thm} \label{2} Let $p$ be a prime number and
$F=\Q(\sqrt{p})$.  Let $H_L$ be the ring class field corresponding
to the order $L=\frak o_F+\frak o_F\sqrt{-1}$. Then the diophantine
equation $x^2+y^2=\alpha$ is solvable over $\frak o_F$ if and only
if $\alpha$ satisfies the Artin condition of $H_L$.
\end{thm}

Now we use Theorem \ref{2} to give an explicit example. Let
$F=\Q(\sqrt{17})$. We write
$N_{F/\Q}(\alpha)=2^{s_1}17^{s_2}p_1^{e_1}\cdots p_g^{e_g}$ for any
$\alpha\in \frak o_F$. Let $D(n)=\{p_1, \cdots, p_g \}$ and
$h(x)=x^4-2x^2+17$. Denote
$$\aligned& D_1=\{p\in D(n): (\frac{-17}{p})=-(\frac{-1}{p})=1\}\cr & D_2=
\{p\in D(n): (\frac{-17}{p})=(\frac{-1}{p})=1 \text{ and } h(x)
\equiv 0 \mod p \text{ is not solvable}\}.
\endaligned $$
We can see $e_i$ is even if $p_i\in D_1$.
\begin{exa} Let $F=\Q(\sqrt{17})$ and let $\alpha$ be an integer in $F$ with
the above notation. Then $x^2+y^2=\alpha$ is solvable over $\frak
o_F$ if and only if

(1) The equation $x^2+y^2=\alpha$ has integral solutions at every
place of $F$.

(2) The sum $$s_1+\sum_{p_i\in D_1}e_i/2+\sum_{p_i\in D_2}e_i\equiv
0 \mod 2.$$
\end{exa}

\bf{Acknowledgment} \it{ The author would like to thank Fei Xu and
Chungang Ji for some helpful discussions. The work is supported by
the Morningside Center of Mathematics and  NSFC, grant \# 10901150
and \# 10671104.}


\end{document}